\documentclass[12pt,twoside]{article}
\usepackage{amsmath}
\usepackage{amsfonts}
\usepackage{amsthm}
\usepackage{amssymb}
\usepackage{graphicx}
\usepackage{multicol}

\begin{document}
\title{{\normalsize{\bf Smooth Homotopy  4-Sphere}}}
\author{{\footnotesize Akio KAWAUCHI}\\
{\footnotesize{\it Osaka Central Advanced Mathematical Institute, Osaka Metropolitan University}}\\
{\footnotesize{\it Sugimoto, Sumiyoshi-ku, Osaka 558-8585, Japan}}\\
{\footnotesize{\it kawauchi@omu.ac.jp}}}
\date\, 
\maketitle
\vspace{0.25in}
\baselineskip=10pt
\newtheorem{Theorem}{Theorem}[section]
\newtheorem{Conjecture}[Theorem]{Conjecture}
\newtheorem{Lemma}[Theorem]{Lemma}
\newtheorem{Sublemma}[Theorem]{Sublemma}
\newtheorem{Proposition}[Theorem]{Proposition}
\newtheorem{Corollary}[Theorem]{Corollary}
\newtheorem{Claim}[Theorem]{Claim}
\newtheorem{Definition}[Theorem]{Definition}
\newtheorem{Example}[Theorem]{Example}

\begin{abstract} Every smooth homotopy  4-sphere is diffeomorphic 
to the  4-sphere. 

\phantom{x}


\end{abstract}

\baselineskip=15pt

\bigskip

\noindent{\bf 1. Introduction}

This paper is a paper created by adding proofs to the research announcement paper \cite{KA}  without proof. So, the proofs are original in this paper. Unless otherwise stated, {\it manifolds}, {\it embeddings}  and  {\it isotopies} are considered in the {\it smooth category}. An $n$-punctured manifold  of  an $m$-manifold  $X$  is the $m$-manifold  $X^{n(0)}$ obtained from  $X$ by removing  the  interiors of  $n$ mutually disjoint $m$-balls in the  interior  of  $X$, where choices of the  $m$-balls are independent of the diffeomorphism type  of  $X^{n(0)}$. 
The  $m$-manifold  $X^{1(0)}$ is denoted by  $X^{(0)}$. 
By this convention, a {\it homotopy  4-sphere} is  a 4-manifold $M$ homotopy equivalent to  
the 4-sphere $S^4$, and a {\it homotopy 4-ball} is a 1-punctured manifold 
$M^{(0)}$ of a homotopy  4-sphere $M$. 
Ever since the positive solution of Topological  4D  Poincar{\'e} Conjecture 
(meaning that every topological homotopy 4-sphere 
is homeomorphic to $S^4$) and existence of exotic 4-spaces, \cite{Fr, FQ,Kirby},  it has 
been questioned whether Smooth 4D  Poincar{\'e} 
Conjecture (meaning that every homotopy 4-sphere 
is diffeomorphic to  $S^4$) holds. This 
paper answers this question affirmatively (see 
Corollary 1.3). 
For a positive integer $n$, the {\it stable 4-sphere of genus} $n$ is  
the  4-manifold
\[\Sigma=\Sigma(n)=S^4\#n(S^2\times S^2)=S^4\#_{i=1}^n S^2\times S^2_i,\]
which is the union of the $n$-punctured manifold $(S^4)^{n(0)}$ of 
the 4-sphere $S^4$ and  the $1$-punctured manifolds 
$(S^2\times S^2_i)^{(0)}\, (i=1,2,\dots, n)$ of the 2-sphere products  
$S^2\times S^2_i\, (i=1,2,\dots, n)$  
pasting the boundary 3-spheres of $(S^4)^{n(0)}$ to the  boundary 3-spheres of 
$(S^2\times S^2_i)^{(0)}\, (i=1,2,\dots, n)$.   
An {\it orthogonal 2-sphere pair} or simply an {\it O2-sphere pair} 
of the stable 4-sphere 
$\Sigma$ is a pair $(S,S')$ of  2-spheres $S$ and $S'$ embedded  in 
$\Sigma$ which meet transversely at  just one point with the intersection numbers 
$\mbox{I}(S,S)=\mbox{I}(S',S')=0$ and $\mbox{I}(S,S')=1$. 
A {\it pseudo-O2-sphere basis} of the stable 4-sphere $\Sigma$ of genus $n$ is 
the system $(S_*,S'_*)$ 
of $n$  disjoint  O2-sphere pairs 
$(S_i,S_i')\, (i=1,2,\dots, n)$ in $\Sigma$.  
Let $N(S_i,S_i')$ be a regular neighborhood of the union $S_i\cup S_i'$ of the 
O2-sphere pair $(S_i,S_i')$ in $\Sigma$ such that $N(S_i,S_i')\, (i=1,2,\dots, n)$ are mutually disjoint and diffeomorphic to the 
1-punctured manifold $(S^2\times S^2)^{(0)}$ 
of the sphere product $S^2\times S^2$.  
The {\it region} of a pseudo-O2-sphere basis 
$(S_*,S'_*)$ in $\Sigma$ of genus $n$ 
is a connected 4-manifold $\Omega(S_*,S'_*)$ in $\Sigma$ obtained from the 4-manifolds 
$N(S_i,S_i')\, (i=1,2,\dots, n)$ by connecting along   disjoint 1-handles 
$h^1_j\, (j=1,2,\dots, n-1)$ in $\Sigma$. 
Since $\Sigma$ is a simply connected 4-manifold, 
the region $\Omega(S_*,S'_*)$ in $\Sigma$  does not depend on any choices of 
the 1-handles  and is uniquely determined
by the pseudo-O2-sphere basis $(S_*,S'_*)$ up to isotopies of $\Sigma$ (see \cite{HK}). 
The {\it residual region} of  $\Omega(S_*,S'_*)$ in $\Sigma$ is the 4-manifold 
${\Omega}^c(S_*,S'_*)=\mbox{cl}(\Sigma\setminus  \Omega(S_*,S'_*))$ which 
is  a  homotopy 4-ball  shown   by  van Kampen theorem  and  a homological  argument.
An {\it O2-sphere basis} of the stable 4-sphere $\Sigma$ of genus $n$ is a 
pseudo-O2-sphere basis $(S_*,S'_*)$ of $\Sigma$ such that the residual region 
${\Omega}^c(S_*,S'_*)$ is diffeomorphic to the 4-ball.  
 For example, the pseudo-O2-sphere basis  
\[(S^2\times1_*,1\times S^2_*) =\{(S^2\times1_i,1\times S^2_i) |\, i=1,2,...,n\}\]
of $\Sigma$ is an O2-sphere basis of $\Sigma$ and called the {\it standard} O2-sphere basis of 
$\Sigma$.
The following result is basically the main result of this paper.

\phantom{x}

\noindent{\bf Theorem~1.1.} For any two pseudo-O2-sphere bases $(R_*,R'_*)$ and 
$(S_*,S'_*)$ of the stable 4-sphere 
$\Sigma$ of any genus $n\geq 1$, there is an orientation-preserving diffeomorphism 
$h:\Sigma\to \Sigma$ sending  $(R_i,R'_i)$ to $(S_i,S'_i)$ for all 
$i\, (i=1,2,\dots,n)$. 

\phantom{x} 

Since the image of an O2-sphere basis $(R_*,R'_*)$ of  $\Sigma$ by  an orientation-preserving 
diffeomorphism $f:\Sigma\to \Sigma$ is  an O2-sphere basis of  $\Sigma$,  
the following corollary is obtained from  existence of the standard O2-sphere basis and Theorem~1.1.

\phantom{x}

\noindent{\bf Corollary~1.2.}  Every  pseudo-O2-sphere basis  of the stable 4-sphere 
$\Sigma$ of any genus $n\geq 1$ is an O2-sphere basis  of $\Sigma$.

\phantom{x}

The following result (4D Smooth Poincar{\'e} Conjecture) 
is a direct consequence of Corollary~1.2. 

\phantom{x}

\noindent{\bf Corollary~1.3.}  
Every homotopy  4-sphere $M$ is diffeomorphic to  the  4-sphere $S^4$.

\phantom{x}

\noindent{\bf Proof of Corollary~1.3.} 
It is known that there is an orientation-preserving diffeomorphism 
\[\kappa: M\# \Sigma\to \Sigma\]
from the connected sum $M\# \Sigma$ onto $\Sigma$ for a positive integer 
$n$ by Wall \cite{W}. 
The connected sum  $M\# \Sigma$ is taken the union $M^{(0)}\cup \Sigma^{(0)}$. 
By Corollary~1.2,  the image $\kappa(\Sigma^{(0)})$ is the region 
$\Omega(S_*,S'_*)$ of an O2-sphere basis $(S_*,S'_*)$ of $\Sigma$.  
and the residual region  
${\Omega}^c(S_*,S'_*)=\kappa(M^{(0)})$ is a 4-ball and hence $M^{(0)}$  is 
diffeomorphic to the 4-ball $D^4$.
The diffeomorphism $M^{(0)}\to D^4$ extends to a diffeomorphism $M\to S^4$ by 
 $\Gamma_4=0$ in Cerf \cite{Cerf} or  $\pi_0(\mbox{Diff}^+(S^3))=0$ by Hatcher \cite{Hat}. 
This completes  the  proof  of   Corollary~1.3.
$\square$

\phantom{x}

As  it  is  seen  from  the proof  of   Corollary~1.3, Theorem~1.1 is actually equivalent to the Smooth 4D  Poincar{\'e} Conjecture by $\Gamma_4=0$ or $\pi_0(\mbox{Diff}^+(S^3))=0$.
It remains unknown whether the diffeomorphism $M \to S^4$ in Corollary~1.3 made orientation-preservingly is isotopically unique, although it is concordantly unique since $\Gamma_5=0$  by Kervaire \cite{Ke}, and piecewise-linear-isotopically unique  by Hudson \cite{Hud}, Rourke-Sanderson \cite{RS}. The Piecewise-Linear 4D Poincar{\'e} Conjecture is equivalent the Smooth 4D Poincar{\'e} Conjecture by using a compatible smoothability of a piecewise-linear 4-manifold and a basic fact that every piecewise-linear auto-homeomorphism of the 4-disk keeping the boundary identically is piecewise-linearly $\partial$-relatively isotopic to the  identity by  \cite{Hud,RS}. The result of Wall \cite{W} used for the proof of Corollary~1.3 says further that for every closed simply connected signature-zero spin 4-manifold $M$ with the second Betti number $\beta_2(M;Z)=2m$, there is a diffeomorphism $\kappa:M\#\Sigma(n)\to \Sigma(m+n)$
for some $n$. Then there is also a homeomorphism  from $M$ to $\Sigma(m)$ by \cite{Fr, FQ}. It should be noted that the present technique used for the proof of Theorem~1.1 cannot be directly generalized to the case of $m>0$. 
In fact, it is known by Akhmedov-Park \cite{AP} that  there is a  closed simply connected 
signature-zero spin 4-manifold $M$ with a large second Betti number 
$\beta_2(M;Z)=2m$ such that $M$ is not diffeomorphic to 
$\Sigma(m)$ . The following corollary is  what can be said in this paper.

\phantom{x}

\noindent{\bf Corollary~1.4.}  Let $M$ and $M'$ be 
any closed (not necessarily simply connected) 
4-manifolds with the same second Betti number 
$\beta_2(M;Z)=\beta_2(M';Z)$. 
Then  an embedding $u: M^{(0)} \to  M'$ induces  the 
fundamental group  isomorphism 
$u_{\#}:\pi_1(M^{(0)},x)\to \pi_1(M',u(x))$
if and only if the embedding $u: M^{(0)} \to M'$ 
extends to a diffeomorphism $u^+: M \to M'$. 

\phantom{x}

\noindent{\bf Proof of Corollary~1.4.}  Since the proof of  the \lq\lq if\rq\rq part
is obvious, it suffices to prove the \lq\lq only if\rq\rq part. 
For this proof, confirm by van Kampen theorem and a homological argument that the closed complement $\mbox{cl}(M'\setminus u(M^{(0)}))$ is a homotopy 4-ball, which is diffeomorphic to 
the 4-ball $D^4$ and  the embedding $u: M^{(0)} \to M'$  extends to a diffeomorphism 
$u^+: M \to M'$ by the proof of Corollary 1.3. This completes the proof of Corollary~1.4.
$\square$
 
\phantom{x}

The following corollary is obtained by combining Corollary~1.3 with 
the triviality condition of an $S^2$-link in $S^4$,  \cite{K'}.

\phantom{x}

\noindent{\bf Corollary~1.5.} Every closed 4-manifold $M$ such that  
the fundamental group $\pi_1(M,x)$ is a free group of rank $n$ 
and $H_2(M;Z)=0$ is diffeomorphic to theclosed 4D handlebody
$Y^S=S^4\#_{i=1}^n  S^1\times S^3_i$.

\phantom{x}

\noindent{\bf Proof of Corollary~1.5.} Let $k_i\, (i=1,2,\dots,n)$ 
be a system of  mutually disjoint simple loops 
in $M$ which is homotopic to a system of loops with leggs to the base point $x$ 
generating the free group $\pi_1(M,x)$, and $N(k_i)= S^1\times D^3_i\, (i=1,2,\dots,n)$ 
a system of mutually disjoint regular neighborhoods of $k_i\, (i=1,2,\dots,n)$ 
in $M$. The 4-manifold $X$ obtained from $M$ by replacing $S^1\times D^3_i$ 
with $D^2\times S^2_i$ for every $i$ is a homotopy 4-sphere by 
van Kampen theorem and $H_2(M;Z)=0$, and hence $X$ is diffeomorphic 
to $S^4$ by Corollary~1.3.  The $S^2$-link 
$L=\cup_{i=1}^n K_i$ in $X=S^4$ with  component  $K_i=S^2_i$  the core 2-sphere of $D^2\times S^2_i$ has the free fundamental group $\pi_1(S^4\setminus L, x)$ of rank $n$ 
with a meridian  basis  since it is canonically isomorphic to 
the free fundamental group  $\pi_1(M,x)$ by a general position argument. 
The $S^2$-link $L$ is a trivial $S^2$-link in $S^4$ and hence 
bounds mutually disjoint 3-balls in $S^4$ by  \cite{K'}.   
By returning $D^2\times S^2_i$ to $S^1\times D^3_i$ for every $i$, 
the 4-manifold $M$ is seen  to be diffeomorphic to the
closed 4D handlebody $Y^S$.  
This completes  the   proof  of    Corollary 1.5.  $\square$

\phantom{x}

The  following corollary (4D Smooth Schoenflies  Conjecture) is also obtained.

\phantom{x}

\noindent{\bf Corollary~1.6.}  Any (smoothly) embedded 3-sphere $S^3$ in 
the 4-sphere $S^4$ splits $S^4$ into two components of 4-manifolds which are 
both diffeomorphic to the 4-ball. 

\phantom{x}

\phantom{x}

\noindent{\bf Proof of Corollary~1.6.}  The splitting components  
are  homotopy 4-balls by  van Kampen theorem and 
a homological argument,  which are diffeomorphic to the 4-ball 
by the proof of Corollary~1.3. 
This  completes  the  proof  of  Corollary 1.6.

 $\square$

\phantom{x}

The paper is organized as follows: In Section~2, a trivial surface-knot in the 4-sphere $S^4$  is discussed to observe that the stable 4-sphere  $\Sigma$  of genus $n$  is the double branched covering space  $S^4(F)_2$ of $S^4$ branched  along a trivial  surface-knot $F$ of  genus $n$.  
An {\it O2-handle basis} of a trivial surface-knot $F$ in $S^4$  is also  
introduced there to show that the lift 
$(S(D_*),S(D'_*))$ of  the core system $(D_*,D'_*)$ of  any O2-handle basis
 $(D_*\times I,D'_*\times I)$ of $F$ to $S^4(F)_2=\Sigma$ is an O2-sphere basis of  $\Sigma$ (See Corollary~ 2.2). In Section~3, the proof  of  Theorem~1.1 is done. In Section~4, any two homotopic diffeomorphisms of the stable 4-sphere $\Sigma$  are shown to be isotopic if one diffeomorphism allows a deformation by an element of 
$\mbox{Diff}^+(D^4,\mbox{rel}\partial)$
(see Theorem 4.1 and Corollaries 4.2, 4.3 for the details).

\phantom{x}

\noindent{\bf 2. Double branched covering of trivial surface-knot for stable 4-sphere}

A {\it surface-knot} of genus $n$ in the 4-sphere $S^4$ is 
a closed surface $F$ of genus $n$ embedded  in $S^4$. 
Two surface-knots $F$ and $F'$ in $S^4$ are {\it equivalent} if there is an 
orientation-preserving diffeomorphism $f: S^4 \to S^4$ sending $F$ to $F'$ 
orientation-preservingly.           
The map $f$ is called  an {\it equivalence}.   
A {\it trivial } surface-knot of genus $n$ in $S^4$ is a surface-knot $F$ of 
genus $n$ which is the boundary of 
a handlebody of genus $n$ embedded  in $S^4$, 
where a handlebody of genus $n$ means a 3-manifold which is a 3-ball for $n=0$, 
a solid torus for $n=1$ or a boundary-disk sum of $n$ solid tori for $n\geq 2$.  
A {\it surface-link} in $S^4$ is a union of disjoint surface-knots in $S^4$, and 
a {\it trivial surface-link} is a surface-link bounding  disjoint  handlebodies 
in $S^4$.  
A trivial surface-link in $S^4$ is  determined regardless of the embeddings 
and unique up to isotopies, \cite{HK}. 

A {\it symplectic basis} of a closed surface $F$ of genus $n$ is a system $(x_*,x'_*)$ 
of element pairs $(x_j,x'_j)\, (j=1,2,\dots,n)$ of $H_1(F;Z)$  with 
the intersection numbers
\[\mbox{I}(x_j, x_{j'})=\mbox{I}(x'_j, x'_{j'})=\mbox{I}(x_j, x'_{j'})=0\] 
for all  $j,j'$  except that $\mbox{I}(x_j, x'_j)=+1$ for all  $j$. 
Every pair $(x_1,x'_1)$ with $\mbox{I}(x_1, x'_1)=+1$ is extended to a 
symplectic basis $(x_*,x'_*)$ of $F$ by an argument on the intersection form 
\[\mbox{I}:H_1(F;Z)\times H_1(F;Z)\to Z.\]
Further, every symplectic basis 
$(x_*,x'_*)=\{(x_j,x'_j)|\,j=1,2,\dots,n\}$ is realized by a system of oriented 
simple loop  pairs $(e_*,e'_*)=\{(e_j,e'_j)|\,j=1,2,\dots,n\}$ of $F$ with   
$e_j\cap e_{j'}=e'_j\cap e'_{j'}=e_j\cap e'_{j'}=\emptyset$ for all distinct $j,j'$  
and with tranverse intersection $e_j\cap e'_j$ at just one point for all  $j$, which is called a {\it loop basis} of  $F$.
For a surface-knot $F$ in $S^4$, 
an element $x\in H_1(F;Z)$ is said to be {\it spin} if the 
$Z_2$-reduction $[x]_2 \in H_1(F;Z_2)$ of $x$ has 
$\eta([x]_2)=0$ for the $Z_2$-quadratic function 
\[\eta:H_1(F;Z_2)\to Z_2\]
associated with a surface-knot $F$ in $S^4$, which is defined as follows:
For  a simple loop $e$ in $F$ bounding  a surface $D_e$ in 
$S^4$ with $D_e\cap F=e$, the $Z_2$-self-intersection number 
$\mbox{I}(D_e,D_e)\pmod{2}$ with respect to the $F$-framing 
is defined to be the value $\eta([e]_2)$. 
For every surface-knot $F$ in $S^4$, there is a spin basis of $F$ (see \cite{HiK}).  
Every spin pair $(x_1,x'_1)$ in $F$ with $\mbox{I}(x_1, x'_1)=+1$ 
is extended to a spin symplectic basis $(x_*,x'_*)$ of $F$ by  vanishing 
of  Arf invarinat of the $Z_2$-quadratic function 
$\eta:H_1(F;Z_2)\to Z_2$ for every surface-knot $F$ in $S^4$. 
In particular, any spin pair $(x_1,x'_1)$ is realized by a spin loop pair 
$(e_1,e'_1)$ of $F$ extendable to a spin loop basis $(e_*,e'_*)$ of $F$. 

A 2-{\it handle} on a surface-knot $F$ in $S^4$ is a 2-handle 
$D\times  I$ on $F$ embedded  in $S^4$ such that 
\[(D\times  I)\cap  F =(\partial  D)\times  I,\]
 where $I$ denotes a closed interval with $0$ as the center 
and $D\times 0$ is called the {\it core} of the 2-handle $D\times I$ and 
identified with $D$.   For a 2-handle $D\times I$ on $F$ in $S^4$, 
the loop $\partial D$ of the core disk $D$ is a spin loop in $F$ since 
$\eta([\partial D]_2)=0$.  
To save notation, if an embedding $h:D\times I\cup F \to X$ is given from a 2-handle 
$D\times I$ on a surface $F$ to a 4-manifold $X$, 
then the 2-handle image $h(D\times I)$ and the core image $h(D)$ on $h(F)$ are denoted by $hD\times I$ and $hD$, respectively.
An {\it orthogonal 2-handle pair} or simply an {\it O2-handle pair} on 
a surface-knot $F$ in $S^4$ is a pair   $(D\times I, D' \times  I)$ 
of 2-handles  $D\times I$ and $D' \times  I$ on  $F$  which {\it meet orthogonally} 
on $F$,   in other words, which 
meet $F$ onlywith the attaching annuli 
$(\partial D)\times I$ and $(\partial D') \times  I$
so that the loops $\partial D$ and $\partial D'$ meet transversely 
at just one point $q$ and the intersection 
$(\partial D)\times I\cap (\partial D') \times  I$ 
is diffeomorphic to the square $Q=\{q\} \times  I\times  I$ 
 \cite{K}.  For a trivial surface-knot F of  genus n in S4, 
 an O2-{\it handle basis} of  $F$ of genus $n$ in $S^4$ is a 
system $(D_*\times I,D'_*\times I)$ of mutually disjoint O2-handle pairs 
$(D_i\times I, D_i'\times I)\, (i=1,2,\dots, n)$ on $F$ in $S^4$ such that 
the loop system 
$(\partial D_*,\partial D'_*)$ given by $\{ (\partial D_i, \partial D_i')|\, i=1,2,\dots,n\}$ 
forms a spin loop basis of $F$. 
Every trivial surface-knot $F$ in $S^4$ is moved into the boundary of 
a standard handlebody in the equatorial 3-sphere $S^3$ of  $S^4$, where   
a {\it standard}  O2-handle basis  and a {\it standard} loop basis of $F$ are taken.
For any given spin loop basis of a trivial  surface-knot $F$ of genus $n$ in 
$S^4$, there is an O2-handle basis $(D_*\times I,D'_*\times I)$ 
of  $F$ in $S^4$ with the given spin loop basis as 
the loop basis $(\partial D_*,\partial D'_*)$. This is because there is an equivalence 
$f:(S^4,F)\to (S^4,F)$  sending  the standard spin loop basis 
to the given spin loop basis of $F$ and hence there is 
an O2-handle basis of $F$ in $S^4$ with the given 
spin loop basis which is the image of the standard O2-handle basis of $F$ by \cite{Hiro}, 
\cite[(2.5.1), (2.5.2)]{K}.  Further  any O2-handle basis of $F$ in $S^4$ with attaching part fixed is unique up to orientation-preserving diffeomorphisms of $S^4$ keeping $F$ point-wise fixed, \cite{KS}. 

For the double branched covering projection  
$p:S^4(F)_2\to S^4$  
branched along $F$,  the non-trivial covering involution of  $S^4(F)_2$ is denoted by $\alpha$.  
The preimage $p^{-1}(F)$  in $\Sigma$ of $F$ which is the fixed point set of $\alpha$  and diffeomorphic to $F$ is also written by  the same notation as $F$.  
The following result is a standard result.

\phantom{x}

\noindent{\bf Lemma~2.1.}  Let $(D_*\times I, D'_*\times I)$ be a standard O2-handle basis  of a trivial surface-knot $F$ of genus $n$ in $S^4$. 
Then there is an orientation-preserving diffeomorphism $f:S^4(F)_2 \to \Sigma$ 
sending the 2-sphere pair system $(S(D_*),S(D'_*))$ with 
$S(D_i)=D_i\cup\alpha D_i$ and $S(D'_i)=D'_i\cup\alpha D'_i\, (i=1,2,\dots,n)$ to the standard O2-sphere basis $(S^2\times 1_*,1\times S^2_*)$ of the stable 4-sphere  
$\Sigma$ of genus $n$. In particular, the 2-sphere pair system $(S(D_*), S(D'_*))$ is an O2-sphere basis of $\Sigma$. 

\phantom{x}

\noindent{\bf  Proof of Lemma~2.1.} 
Let $A_i\, (i=1,2,\dots,n)$ be mutually disjoint 4-balls which are regular neighborhoods of the 3-balls 
\[D_i\times I\cup D_i'\times I\quad (i=1,2,\dots,n)\] 
in  $S^4$. The closed complement 
$(S^4)^{n(0)}=\mbox{cl}(S^4\setminus \cup_{i=1}^n A_i)$ 
is the $n$-punctured  manifold of $S^4$. 
Let $P= F\cap (S^4)^{n(0)}$ be a proper $n$-punctured 2-sphere in $(S^4)^{n(0)}$. 
Since the pair $((S^4)^{n(0)} , P)$ 
is  an $n$-punctured pair of a trivial 2-knot space $(S^4,S^2)$ and 
the double branched covering space $S^4(S^2)_2$ is diffeomorphic to  $S^4$, 
the double branched covering space $(S^4)^{n(0)}(P)_2$ of  $(S^4)^{n(0)}$  branched along $P$ 
is diffeomorphic to $(S^4)^{n(0)}$. 
On the other hand, for the proper surface $P_i=F\cap A_i$ in the 4-ball $A_i$, 
the pair $(A_i,P_i)$ is considered as a 1-punctured pair of  a trivial torus-knot 
space $(S^4,T)$, 
so that the double branched covering space 
$A_i(P_i)_2$ is diffeomorphic to the 1-punctured manifold of 
the  double branched covering space $S^4(T)_2$. 
 The trivial torus-knot space $(S^4,T)$ is the double of the product pair 
$(B,o)\times I=(B\times I,o\times I)$ for a trivial loop $o$ in the interior of 
a 3-ball $B$ and an interval $I$, so that 
$(S^4,T)$ is diffeomorphic to the boundary pair 
\[\partial((B,o)\times I^2)=(\partial (B\times I^2),\partial (o\times I^2) ),\]
where $I^m$ denotes the $m$-fold product of $I$ for any $m\geq 2$. 
Thus, the double branched  covering  space $S^4(T)_2$  is diffeomorphic to 
the boundary $\partial (B(o)_2\times I^2)$, where $B(o)_2$ is 
the double branched covering space of $B$ branched along $o$ 
which is diffeomorphic to the product $S^2\times I$. 
This means that  the 5-manifold $B(o)_2\times I^2$ is  the product $S^2\times I^3$,   
and  $S^4(T)_2$ is diffeomorphic to $S^2\times S^2$, 
which shows that $A_i(P_i)_2$ is diffeomorphic to $(S^2\times S^2)^{(0)}$. 
This construction also shows that there is an orientation-preserving diffeomorphism 
\[f_i:A_i(P_i)_2\to (S^2\times S^2)^{(0)}_i\] 
sending the O2-sphere pair $(S(D_i),S(D'_i))$ 
to the standard O2-sphere pair $(S^2\times 1_i,1\times S^2_i)$ of 
the connected summand $(S^2\times S^2)^{(0)}_i$ of  $\Sigma$ for all $i$.  
A desired orientation-preserving diffeomorphism  $f : S^4(F)_2\to \Sigma$ is constructed from 
a diffeomorphism $f'': (S^4)^{n(0)}(P)_2 \to  (S^4)^{n(0)}$ and the diffeomorphisms  $f_i\, (i=1,2,\dots,n)$. 
This completes  the   proof of Lemma~2.1.

$\square$

\phantom{x}

The identification of  $S^4(F)_2=\Sigma$ is fixed  by  
an orientation-preserving diffeomorphism  $f:S^4(F)_2\to \Sigma$ 
given in Lemma~2.1. The following corollary is obtained from Lemma~2.1  
and \cite{K,KS}. 

\phantom{x}

\noindent{\bf Corollary~2.2.} For any two O2-handle bases $(D_*\times I, D'_*\times I)$ and 
$(E_*\times I,E'_*\times I)$ of  a  trivial surface-knot $F$  of  genus $n$  in $S^4$, there is an orientation-preserving $\alpha$-equivariant diffeomorphism $f"$ of   $\Sigma$ sending the  
2-sphere pair system $(S(D_*),S(D'_*))$ to the 2-sphere pair system $(S(E_*),S(E'_*))$. 
In particular, the 2-sphere pair system $(S(D_*),S(D'_*))$ for every 
O2-handle basis $(D_*\times I, D'_*\times I)$ is an O2-sphere basis of $\Sigma$. 

\phantom{x}

\noindent{\bf Proof of Corollary~2.2.} 
There is an equivalence $f:(S^4,F)\to (S^4,F)$ keeping $F$ set-wise fixed 
sending the O2-handle basis   $(D_*\times I, D'_*\times I)$ to the O2-handle basis  
$(E_*\times I,E'_*\times I)$   of  $F$ by uniqueness of an O2-handle pair, \cite{K,KS}. 
By construction, the lifting diffeomorphism $f '':S^4(F)_2\to S^4(F)_2$ of  $f$ 
sends $(S(D_*),S(D'_*))$ to $(S(E_*),S(E'_*))$. 
From Lemma~2.1,  the  2-sphere pair system $(S(D_*),S(D'_*))$ for every 
O2-handle basis $(D_*\times I, D'_*\times I)$ is shown  to  be an O2-sphere basis of $\Sigma$ by taking . $(E_*\times I,E'_*\times I)$ a standard O2-handle basis of $F$.
This completes the  proof of Corollary~2.2. 
$\square$

\phantom{x}

An $n$-{\it rooted disk family} is the triplet $(d, d_*, b_*)$ where $d$ is a disk, 
$d_*$ is a system  of  mutually disjoint disks $d_i\,(i=1,2,\dots,n)$ in the interior of $d$ 
and $b_*$ is a system  of  mutually disjoint bands $b_i\,(i=1,2,\dots,n)$ 
in the $n$-punctured disk $\mbox{cl}(d\setminus d_*)$ such that $b_i$ spans an arc in
the loop $\partial d_i$ and an arc in the loop $\partial d$. 
Let  $\gamma( b_*)$ denote the centerline system of the band system $b_*$.
In the  following lemma shows that there is a canonical $n$-rooted disk family 
$(d, d_*, b_*)$ associated with an O2-handle basis $(D_*\times I, D'_*\times I)$ of a 
trivial surface-knot $F$ of genus $n$ in $S^4$.

\phantom{x}

\noindent{\bf Lemma~2.3.} Let 
$(D_*\times I, D'_*\times I)$ be an  O2-handle basis of  a trivial  surface-knot 
$F$ of genus $n$ in $S^4$, and $(d, d_*,b_*)$ an $n$-rooted disk family.
Then there is an embedding 
\[\varphi:(d, d_*, b_*)\times I \to (S^4, D_*\times I, D'_*\times I)\] 
such that 

\medskip

\noindent{(1)} The surface $F$ is the boundary of 
the handlebody $V$ of genus $n$ given by  
\[V=\varphi(\mbox{cl}(d\setminus d_*)\times I),\] 

\medskip

\noindent{(2)} There is an identification 
\[\varphi(d_*\times I, d_*)=(\varphi(d_*)\times I,\varphi(d_*)) =(D_*\times I,D_*)\] 
as  2-handle systems on $F$, and

\medskip

\noindent{(3)} There is an identification 
\[\varphi(b_*\times I, \gamma(b_*) \times I)=(D'_*\times I,D'_*) \] 
as  2-handle systems on $F$. 

\phantom{x}

Lemma~2.3 says that the 2-handle system s$D_*\times I$ and $D'_*\times I$ are 
attached to the handlebody $V$ bounded by $F$ along a longitude system and a meridian system of $V$, respectively.

\phantom{x}

\noindent{\bf Proof of Lemma~2.3.} If $(E_*\times I, E'_*\times I)$ is 
a standard O2-handle basis of $F$, then it is easy to construct such an embedding 
\[\varphi':(d, d_*, b_*)\times I \to (S^4, E_*\times I, E'_*\times I)\] 
with (1)-(3) taking $\varphi'$  and $(E_*\times I, E'_*\times I) $ as $\varphi$ 
and $(D_*\times I, D'_*\times I)$, respectively. 
In general, there is an equivalence $f:(S^4,F)\to (S^4,F)$ keeping $F$ set-wise fixed 
and sending the  standard  O2-handle basis $(E_*\times I, E'_*\times I)$ of  $F$ into 
the O2-handle basis $(D_*\times I, D'_*\times I)$ of $F$ 
by uniqueness of an O2-handle pair in \cite{K,KS}.  
The composite embedding 
\[\varphi=f\varphi':(d, d_*, b_*)\times I \to (S^4, D_*\times I, D'_*\times I)\] 
is a desired embedding. This completes the proof of Lemma~2.3. 
$\square$

\phantom{x}

 In Lemma~2.3, the embedding $\varphi$,  the 3-ball $B=\varphi(D\times I)$, the 
handlebody $V$  and the pair $(B,V)$ are respectively  called a {\it bump embedding}, 
a {\it bump 3-ball}, a {\it bump handlebody} and a {\it bump pair} of $F$ in $S^4$. 
For a bump embedding 
\[\varphi:(d, d_*,b_*)\times I \to (S^4, D_*\times I,D'_*\times I),\] 
there is an embedding 
$\varphi'':d\times I \to S^4(F)_2$ with $p\varphi''=\varphi$. 
Since $p(\varphi''(d_*\times I),\varphi''(b_*\times I))=(D_*\times I,D'_*\times I)$ 
by the conditions (1)-(3) of Lemma~2.3, 
the images $\varphi''(d_*\times I)$ and $\varphi''(b_*\times I)$ 
are respectively considered as 2-handle systems on $F$ in $S^4(F)_2$ with  
$p\varphi''(d_*\times I)=\tilde D_*\times I$ and $p\varphi''(b_*\times I)=\tilde D'_*\times I$ 
so that  $(\varphi''(d_*\times I),\varphi''(b_*\times I))$ is an O2-handle basis 
of $F$ in $S^4(F)_2$, which is also denoted by $(D_*\times I,D'_*\times I)$ 
to define an embedding 
\[\varphi'': (d,d_*,b_*)\times I \to (S^4(F)_2,D_*\times I,D'_*\times II)\]
with $p\varphi''=\varphi$. This embedding is called a {\it lifting bump embedding}  of  the  bump embedding.  In this case, the {\it bump 3-ball} 
 $\varphi''(d\times I)$  and the {\it bump handlebody} 
 $\varphi''(\mbox{cl}(d\setminus d_*)\times I)$  of  $F$ in $S^4(F)_2$  are
also denoted by $B$ and $V$ by counting 
$p\varphi''(d\times I)=B$ and  
 $p\varphi''(\mbox{cl}(d\setminus d_*)\times I)=V$, respectively. 
 For the non-trivial covering involution $\alpha$ of  $S^4(F)_2$,
the composite embedding 
\[\alpha \varphi'' :(d, d_*,b_*)\times I \to 
(S^4(F)_2, \alpha D_*\times I, \alpha D'_*\times I)\]
is another lifting bump embedding of the bump embedding $\varphi$. 
For the bump 3-ball $\alpha \varphi''(d\times I)=\alpha B$ and  
the bump handlebody $\alpha\varphi'' (\mbox{cl}(d\setminus d_*)\times I)=\alpha V$  
of $F$  in $S^4(F)_2$,  we have 
\[ V\cap \alpha V=B\cap \alpha B=F\]
in $S^4(F)_2$. 
For a  lifting bump embedding, the following lemma is obtained.

\phantom{x}

\noindent{\bf  Lemma~2.4.} 
Let  $\varphi'' :   (d, d_*,b_*)\times I \to (\Sigma, D_*\times I,D'_*\times I)$  
be a lifting bump embedding. For an embedding 
\[u:\Sigma^{(0)}\to \Sigma,\]
assume that the image $\varphi''(d\times I)$ is in the interior of $\Sigma^{(0)}$ 
to define the composite embedding $u \varphi'': (d, d_*,b_*)\times I \to 
(\Sigma, uD_*\times I,uD'_*\times I)$.  
Then there is a diffeomorphism $g:\Sigma\to\Sigma$ which is isotopic to the identity 
such that the composite embedding 
\[gu \varphi'': (d, d_*,b_*)\times I \to 
(\Sigma, guD_*\times I,guD'_*\times I)\]
 is identical to the lifting bump embedding  
\[\varphi'': (d, d_*,b_*)\times I \to (\Sigma, D_*\times I,D'_*\times I). \]

\phantom{x}

\noindent{\bf Proof of Lemma~2.4.}  
The 0-section $(d, d_*,b_*)\times 0$ of the line bundle $(d, d_*,b_*)\times I$ 
of the $n$-rooted disk family $(d,d_*,b_*)$  is identified with $(d,d_*,b_*)$. 
Move the disk $u\varphi''(d)$ into the disk $\varphi''(d)$ in $\Sigma$ and then move the disk system $u\varphi''(d_*)$ and the band system $u\varphi''(b_*)$ into the disk system 
$\varphi''(d_*)$ and 
the band system $\varphi''(b_*)$ in the disk $\varphi''(d)$, respectively.
These deformations are attained by a diffeomorphism $g'$ of $\Sigma$ 
which is isotopic to the identity, so that  
 \[g' u \varphi''(d, d_*,b_*)=\varphi''(d, d_*,b_*).\]
Further,  there is a diffeomorphism $g''$ of $\Sigma$ which is isotopic to the identity 
such that the composite embedding $g''g'u:\Sigma^{(0)}\to \Sigma$ preserves 
the normal line bundles of the disk $\varphi''(d)$ in $\Sigma^{(0)}$ and $\Sigma$,  
so that  
\[g u \varphi''(d, d_*,b_*)\times I=\varphi''(d, d_*,b_*)\times I\] 
for the diffeomorphism $g=g''g'$ of $\Sigma$  isotopic to 
the identity. This completes  the proof of Lemma~2.4. $\square$ 

\phantom{x}

In Lemma~2.4, also assume  that the image $\alpha\varphi''(d\times I)$ is in the interior of 
$\Sigma^{(0)}$. Then note that 
any disk interior of the disk systems $gu\alpha D_*$ and 
$gu\alpha D'_*$ does not meet the bump 3-ball $B=\varphi''(d\times I)$ in $\Sigma$. 
In fact,  $gu$ defines an embedding from 
\[gu: B\cup \alpha B  \to \Sigma\]
 with $gu(B,F)=(B,F)$.    
The complement $gu\alpha B \setminus F$ of $F$ in the  3-ball  $gu\alpha B $ does not meet the bump 3-ball $B$ since $B\cap \alpha B  =F$. 
This means that any disk interior  of the disk systems $gu\alpha D_*$ and $gu\alpha D'_*$ does not meet the bump 3-ball $B$. 
Note that this property comes from  the  fact that $\Sigma^{(0)}$ and $\Sigma$ have the same  genus $n$.

\phantom{x}

\noindent{\bf 3. Proof of  Theorem~1.1}

For the proof of Theorem~1.1, the following result known by Wall \cite{W0}  is used.

\phantom{x}

\noindent{\bf Lemma~3.1.}  
For any  pseudo-O2-sphere bases $(R_*,R'_*)$ and  
$(S_*,S'_*)$ of the stable 4-sphere $\Sigma$ of genus $n$, 
there is an  orientation-preserving diffeomorphism 
$f :\Sigma\to \Sigma$ which induces an isomorphism 
$f_*:H_2(\Sigma;Z) \to H_2(\Sigma;Z)$ 
such that  
$[f R_i]=[S_i]$ and $[f R_i']=[S'_i]$
for all $i$. 

\phantom{x}

Assum that $(R_*,R'_*)$ is an O2-sphere basis of  $\Sigma$ with $(R_*,R'_*)=(S(D_*),S(D'_*))$   
for an O2-handle basis  
$(D_*\times I, D'_*\times I)$ of  a trivial  surface-knot $F$ of genus $n$ in $S^4$.  Let 
$u:\Sigma^{(0)}\to \Sigma$ be an embedding such that  $(uS(D_*),uS(D'_*))=(S_*,S'_*)$. 
By Lemma~3.1, assume that the homology classes $[uS(D_i)]=[S_i]$  and 
$ [uS(D'_i)]=[S'_i]$  are respectively identical to the homology classes 
$[Ri]=[S(Di)]$ and $[R'i]=[S(D'i)]$ for all $i$. 
 Let  $(B,V)$ be a bump pair of the O2-handle basis  
$(D_*\times I, D'_*\times I)$ of $F$ in $S^4$ defined soon after Lemma~2.3. 
Recall that the two lifts  of $(B,V)$ to $\Sigma$ under the double branched covering projection $p:S^4(F)_2\to S^4$ are denoted by $(B,V)$ and $(\alpha B ,\alpha V)$.  
To complete the proof of Theorem~1.1, three lemmas are provided from here. The first lemma is stated as follows.

\phantom{x}

\noindent{\bf Lemma~3.2.}  There is a  diffeomorphism $g$ of 
$\Sigma$ which is isotopic to the identity 
such that the composite embedding 
\[g u:\Sigma^{(0)}\to \Sigma\]
preserves the bump pair $(B,V)$ in $\Sigma$  identically 
and has the property that  every disk interior of the disk systems 
$g u\alpha  D_*$ and $g u\alpha  D'_*$
meets  every disk of the disk systems $\alpha  D_*$ and $\alpha  D'_*$ transversely 
only  with the intersection number $0$.

\phantom{x}

\noindent{\bf Proof of Lemma~3.2.}
By Lemma~2.4, there is a diffeomorphism $g:\Sigma\to\Sigma$ 
which is isotopic to the identity 
such that the composite embedding  $g u:\Sigma^{(0)}\to \Sigma$
preserves the bump pair $(B,V)$ in $\Sigma^{(0)}$ 
identically and has the property that 
any disk interior of the disk systems $gu\alpha D_*$ and $gu\alpha D'_*$
does not meet the O2-handle basis $(D_*\times I,D'_*\times I)$ in $\Sigma$ 
and meets transversely any disk interior of the disk systems $\alpha D_*$ and $\alpha D'_*$ 
with a finite number of points. Since 
\[S(D_i)=D_i\cup\alpha D_i,\,  S(D'_i)=D'_i\cup\alpha D'_i,\,guD_i=D_i\,, guD'_i=D'_i\]
and any disk interior pair of the disk systems $\alpha D_*$ and $\alpha D'_*$ is a disjoint pair, 
every disk interior of  the  disk systems  $gu\alpha D_*$ and $gu\alpha D'_*$ meets every disk interior of the disk systems $\alpha D_*$ and $\alpha D'_*$ only with  intersection number $0$ 
by the homological identities 
\[[guS(D_i)]=[S(D_i)], \quad [guS(D_i)]=[S(D_i)]\]
for all $i$ and the invariance of their intersection numbers. 
This completes the proof of Lemma~3.2. 
$\square$

\phantom{x}

By Lemma~3.2,  assume that the  orientation-preserving embedding 
$u:\Sigma^{(0)}\to \Sigma$ 
sends the bump pair $(B,V)$ to itself identically and has the property that 
every disk interior of the disk systems $u\alpha D_*$ and 
$u\alpha D'_*$ meets  every disk interior  of  the disk systems $\alpha  D_*$ and $\alpha  D'_*$ only with  the intersection number $0$. Then the following lemma is obtained: 

\phantom{x}

\noindent{\bf Lemma~3.3.} There is a diffeomorphism $g$ of $\Sigma$ 
which is isotopic to the identity such that 
the composite embedding $g u:\Sigma^{(0)} \to \Sigma$ 
sends the disk systems $D_*$ and $D'_*$ identically 
and the disk interiors of the disk systems $g u\alpha D_*$,  
$g u\alpha D'_*$ to be  disjoint from the disk systems 
$\alpha  D_*$ and $\alpha  D'_*$ in $\Sigma$. 

\phantom{x}

\noindent{\bf Proof of Lemma~3.3.}
Between the open disks $u\mbox{Int}(\alpha D_i)$, 
$u\mbox{Int}(\alpha D'_{i'})$ for all $i,i'$ and the open disks 
$\mbox{Int}(\alpha D_j)$,  $\mbox{Int}(\alpha D'_{j'})$ for all $j,j'$, 
suppose  an  open disk, say $u\mbox{Int}(\alpha D_i)$ 
meets  an open disk, say $\mbox{Int}(\alpha D_j)$ with a pair of points 
with opposite signs. A procedure to eliminate this pair of points is explained 
from now. 
Let $a$ be a simple arc in the open disk $\mbox{Int}(\alpha D_j)$ 
joining the pair of points whose interior does not meet the intersection points
$\mbox{Int}(\alpha D_j)\cap u\mbox{Int}(\alpha D_*)$ and 
$\mbox{Int}(\alpha D_j)\cap u\mbox{Int}(\alpha D'_*)$.
Let $T(a)$ be the torus obtained from the 2-sphere $uS(D_i)=D_i\cup u(\alpha D'_i)$
by a surgery along a 1-handle $h^1(a)$  on $uS(D_i)$ with core the arc $a$ and with 
$h^1(a) \cap \mbox{Int}(\alpha D_j)=a$. 
Slide the arc $a$  along  the open disk $\mbox{Int}(\alpha D_j)$  
without moving the endpoints so that 

\medskip

\noindent{(*)} the 1-handle $h^1(a)$  passes once time through 
a thickening $S'_j\times I$ of  a 2-sphere $S'_j$ 
parallel to the 2-sphere $S(D'_j)=D'_j\cup \alpha D'_j$ (not meeting  $S(D'_j)$), and 

\medskip

\noindent{(**)}  the intersection $h^1(a)\cap (S'_j\times I)$ is  a 2-handle 
$d_j\times I$ on the torus $T(a)$ which is  a strong deformation 
retract of the 2-handle $h^1(a)$ on $T(a)$. 

\phantom{x}

After these deformations (*), (**), let $u'S(D_i)$ be the 2-sphere obtained  from $T(a)$ by the surgery along the 
2-handle $d_j\times I$. The resulting 2-sphere  $u'S(D_i)$ is obtained from 
the 2-sphere $uS(D_i)$ as $u'S(D_i)=g' u S(D_i)$ for
a diffeomorphism $g':\Sigma\to \Sigma$ which is 
isotopic to the identity.  
This isotopy of $g'$ keeps the outside of a regular neighborhood 
of $h^1(a)$ in the image of $g'$ fixed.  
Next, regard the 2-handle $d_j\times I$  on $T(a)$ as a 1-handle of the 
2-sphere $u'S(D_i)$. 
Let $u''S(D_i)$ be the 2-sphere obtained from $T(a)$  
by the surgery along the 2-handle $\mbox{cl}(S'_j\setminus d_j)\times I$ 
and regard this 2-handle as a 1-handle on the 2-sphere $u''S(D_i)$. 
The 2-sphere $u'S(D_i)$ is isotopically deformed  into the 2-sphere $u''S(D_i)$ by a homotopy deformation  of   a 1-handle (see \cite[Lemma~1.4]{HK}).  
Thus, there is  a diffeomorphism $g'':\Sigma\to \Sigma$ 
which is isotopic to the identity such that  
\[u''S(D_i)=g'''u' S(D_i)=g'' g' uS(D_i).\]
This isotopy of $g''$ keeps the outside of a regular neighborhood of 
$S'_j\times I$ fixed. 
By this procedure, the total geometric intersection number between 
the open disks 
$u''\mbox{Int}(\alpha D_i)$, $u''\mbox{Int}(\alpha D'_{i'})$ for all $i,i'$ 
and the open disks  $\mbox{Int}(\alpha D_j)$,  
$\mbox{Int}(\alpha D'_{j'})$ for all $j,j'$ is reduced by 2. 
By continuing this process, 
we have a diffeomorphism $g:\Sigma\to \Sigma$ 
which is  isotopic to the identity such that 
the composite embedding  $g u:\Sigma^{(0)} \to \Sigma$ sends  
the disk systems $D_*$ and $D'_*$ identically   and the open disks 
$gu\mbox{Int}(\alpha D_i)$ and   
$gu\mbox{Int}(\alpha D'_{i'})$ are disjoint from  the open disks 
$\mbox{Int}(\alpha D_j)$ and $\mbox{Int}(\alpha D'_{j'})$ for all $i,i',j,j'$. 
The procedure is similarly done for the other cases 
that $u\mbox{Int}(\alpha D_i)$ meets  $\mbox{Int}(\alpha D'_{j'})$ 
with a pair of points with opposite signs, 
that  $u\mbox{Int}(\alpha D'_{i'})$ meets  $\mbox{Int}(\alpha D_j)$ 
with a pair of points with opposite signs  and that
$u\mbox{Int}(\alpha D'_{i'})$ meets  $\mbox{Int}(\alpha D'_{j'})$ with a pair of points 
with opposite signs. 
This completes the proof of Lemma~3.3. $\square$ 

\phantom{x}

For the O2-sphere basis $(S(D_*),S(D'_*))$ of $\Sigma$, let 
\[q_*=\{q_i=S(D_i)\cap S(D'_i)|\, i=1,2,\dots,n\}\] 
be the transverse intersection point system  between $S(D_*)$ and $S(D'_*)$.  
The diffeomorphism $g$ of $\Sigma$ sending the disk systems $D_*$ and 
$D'_*$ identically in Lemma~3.3 is  further deformed so that, while leaving the transverse intersection point $q_i$,   the disks $guD_i$ and $D_i$ are separated and then  the disks 
$guD'_i$ and $D'_i$ are separated.  Thus, 
\[guD_i\cap D_i=guD'_i\cap D'_i=q_i \] 
for all $i$. 
By this deformation, the pseudo-O2-sphere basis $(guS(D_*),guS(D'_*))$ of $\Sigma$
 is assumed to meet the O2-sphere basis $(S(D_*),S(D'_*))$ of $\Sigma$ at just the transverse intersection point system $q_*$.  
Next, the diffeomorphism $g$ of $\Sigma$ is deformed  so that 
a disk neighborhood system of $q_*$ in $g uS(D_*)$  and 
a disk neighborhood system of $q_*$ in $S(D_*)$ are matched, and then  
a disk neighborhood system of 
$q_*$ in $g uS(D'_*)$  and a disk neighborhood system of $q_*$ in $S(D'_*)$ 
are matched.
Thus, there is a diffeomorphism $g$ of $\Sigma$ which is isotopic to the identity 
such that the meeting part of the pseudo-O2-sphere basis  
$(g uS(D_*),g uS(D'_*))$ and the O2-sphere basis $(S(D_*),S(D'_*))$ 
is just a disk neighborhood pair system $(d_*,d'_*)$ of the 
transverse intersection point system $q_*$.
Now, assume that for an embedding $u:\Sigma^{(0)} \to \Sigma$,  
the meeting part of the pseudo-O2-sphere basis  
$(uS(D_*),uS(D'_*))$ and the O2-sphere basis $(S(D_*),S(D'_*))$ 
is just a disk neighborhood pair system $(d_*,d'_*)$ of $q_*$.
Then the following lemma is obtained:

\phantom{x}

\noindent{\bf Lemma~3.4.} There is  an orientation-preserving  diffeomorphism 
$h$ of $\Sigma$  such that the composite embedding  
$h u:\Sigma^{(0)} \to \Sigma$ 
preserves the O2-sphere basis $(S(D_*),S(D'_*))$ identically. 

\phantom{x}

The proof of Lemma~3.4 is obtained by using Lemma~3.5 
( {\it Framed Light-bulb Diffeomorphism Lemma})  which is proved easily 
in comparison with an isotopy version (Lemma~3.7)  of this lemma using 
Gabai's 4D light-bulb theorem \cite{G}.  To state Lemmas~3.5,3.7, call 
  a 4-manifold $Y$ in $S^4$ which is diffeomorphic to 
$S^1\times D^3$ a {\it 4D solid torus}. A {\it boundary fiber circle} of the 4D solid torus $Y$ 
is a fiber circle of the $S^1$-bundle $\partial Y$ diffeomorphic to $S^1\times S^2$. 
Let $Y^c=\mbox{cl}(S^4\setminus Y)$ be the exterior of $Y$ in $S^4$.
Let $Y_*$ be a system of mutually disjoint 4D solid tori $Y_i,\, (i=1,2,\dots, n)$ 
in $S^4$, and $Y^c_*$ the system of the exteriors $Y^c_i$ of $Y_i$ in $S^4\, (i=1,2,\dots, n)$.  
Let 
\[\cap Y^c_*=\cap_{i=1}^n Y^c_i.\]

\phantom{x}

\noindent{\bf Lemma~3.5 (Framed Light-bulb Diffeomorphism Lemma).}  
Let $Y_*$ be a system  of  mutually disjoint 4D solid  tori $Y_i\, (i=1,2,\dots,n)$  in $S^4$. 
Let $D_*\times I$ and $E_*\times I$ be  systems of   mutually disjoint framed disks 
 $D_i\times I, E_i\times I\, (i=1,2,\dots,n)$ in $\cap Y^c_*$ such that 
\[(D_*\times I)\cap \partial Y^c_i =(\partial D_i)\times I= (\partial E_i)\times I 
= (E_*\times I)\cap \partial Y^c_i\] 
and $\partial D_i =\partial E_i$  is a  boundary fiber circle  of $Y_i$ for all  $i$. 
Then there is an orientation-preserving diffeomorphism $h:S^4 \to S^4$ sending $Y_*$ identically such that $h(D_*\times I,D_*)=(E_*\times I,E_*)$.

\phantom{x}

\noindent{\bf Proof of Lemma~3.5.}
Let $k_*$ be the system of the loops 
$k_i=\partial D_i=\partial E_i\, (i=1,2,\dots,n)$.
Let $c:k_*\times[0,1]\to D_*$ be a boundary collar function of $D_*$
with $c(x,0)=x$ for all $x\in k_*$, and 
$c':k_*\times[0,1]\to E_*$ a boundary collar function of $E_*$ 
with $c'(x,0)=x$ for all $x\in k_*$.
Assume that 
\[c(x,t)\times I=c'(x,t)\times I\] 
for all $x\in k_*$ and $t\in [0,1]$. 
Let 
\[\nu(\partial D_*)=c(k_*\times[0,1])\quad \mbox{and}\quad  
D^-_*=\mbox{cl}(D_*\setminus \nu(\partial D_*)).\] 
Similarly, let 
\[\nu(\partial E_*)=c(k_*\times[0,1])\quad \mbox{and}\quad  
E^-_*=\mbox{cl}(E_*\setminus \nu(\partial E_*)).\] 
Consider  $D_*\times I$ and $E_*\times I$ in $S^4$. 
Let $\beta_*$ be the system of arcs $\beta_i\, (i=1,2,\dots,n)$ such that 
$\beta_i$ is an arc in $k_i$, 
and $\beta_*^c$ the system of the arcs 
$\beta^c_i=\mbox{cl}(k_i\setminus \beta_i)\, (i=1,2,\dots,n)$.
By 3-cell moves within a regular neighborhood of 
$c(\beta_*\times[0,1])\times I \cup D^-_*\times I$ in $S^4$,  
there is an orientation-preserving diffeomorphism $h'$ of $S^4$ such that 
$h'(D_*\times I)= c(\beta^c_*\times[0,1])\times I$. 
Since $\nu(\partial D_*)$ is identical to $\nu(\partial E_*)$,  
$c(\beta^c_*\times[0,1])$ is identical to $c'(\beta^c_*\times[0,1])$. 
By  3-cell moves within a regular neighborhood of 
$c(\beta_*\times [0,1])\times I \cup E^-_*\times I$ in $S^4$,  
there is an orientation-preserving diffeomorphism $h''$ of $S^4$ such that
$h''(c'(\beta^c_*\times[0,1])\times I)= E_*\times I$. 
The diffeomorphism $h''h'$ is an orientation-preserving diffeomorphism
of $S^4$ sending  $D_*\times I$ to  $E_*\times I$. 
Let $N(k_*)$ be a regular neighborhood system of the loop system 
$k_*$ in $S^4$ meeting $c(k_*\times[0,1])$ regularly,  which is a system of 
$n$ mutually disjoint 4D solid tori. 
Then the diffeomorphism $h''h'$ is deformed into an orientation-preserving 
diffeomorphism $h$ of $S^4$ which sends  $N(k_*)$ identically such that 
$h(D_*\times I,D_*)=(E_*\times I,E_*)$. Here, the 4D solid torus system $N(k_*)$
can be replaced with any given $Y_*$ because $N(k_*)$ is isotopic to $Y_*$ in $S^4$. 
This completes  the proof of Lemma~3.5. $\square$ 

\phantom{x}
 
The proof of Lemma~3.4 is obtained from Lemma~3.5 as follows.

\phantom{x}

\noindent{\bf Proof of Lemma~3.4.}
Consider the 4-manifold $X$ diffeomorphic to  the 4-sphere $S^4$ which is obtained from 
$\Sigma$ by replacing a regular neighborhood system 
$N(S(D'_*))=S^2\times D^2_*$ of the 2-sphere system $S(D'_*)$ in $\Sigma$ with 
the 4D solid torus system $Y_*=D^3\times S^1_*$.  
Let $E^u_*\times I$ and $E_*\times I$ be the 2-handle systems in $X=S^4$ 
attached to the 4D solid torus system $Y_*$ which are
obtained from the thickening 2-sphere systems $uS(D_*)\times I$ and 
$S(D_*)\times I$ in $\Sigma$, respectively. 
Lemma~3.5 can be used for the 2-handle systems $E^u_*\times I$ and 
$E_*\times I$ attached to $Y_*$. Then, 
there is an orientation-preserving  diffeomorphism $\rho:S^4\to S^4$ sending 
$Y_*$ identically  such that  
\[\rho (E^u_*\times I,E^u_*)=(E_*\times I,E_*).\] 
By returning the 4D solid torus system $Y_*$ in $X$ to the regular neighborhood system 
$N(S(D'_*))$ of the 2-sphere system $S(D'_*)$ in $\Sigma$, 
there is an orientation-preserving  diffeomorphism 
$\rho':\Sigma\to \Sigma$ sending $N(S(D'_*))$ identically  such that 
$\rho' u S(D_*)=S(D_*)$.  
For the pseudo-O2-sphere basis  $(S(D_*), \rho' u S(D'_*))$
and the O2-sphere basis $(S(D_*),S(D'_*))$ in $\Sigma$, consider 
the 4-manifold $X'$ obtained from $\Sigma$ by replacing a regular neighborhood system 
$N(S(D_*))=S^2\times D^2_*$ of the 2-sphere system $S(D_*)$ in $\Sigma$ with 
the 4D solid torus system $Y'_*=D^3\times S^1_*$. 
Then the 4-manifold $X'$ is diffeomorphic to the 4-sphere $S^4$.
Let $E'{}^u_*\times I$ and $E'_*\times I$ be the 2-handle systems in $X'=S^4$ 
attached to the 4D solid torus system $Y'_*$ which are
obtained from the thickening 2-sphere systems $\rho' uS(D'_*)\times I$ and 
$S(D'_*)\times I$ in $\Sigma$, respectively.     
Lemma~3.5 can be used for the 2-handle systems $E'{}^u_*\times I$ and $E'_*\times I$ 
attached to $Y'_*$. Then, 
there is an orientation-preserving  diffeomorphism $\rho'':S^4\to S^4$ sending 
$Y'_*$ identically  such that  
\[\rho'' (E'{}^u_*\times I,E'{}^u_*)=(E'_*\times I,E'_*).\] 
By returning the 4D solid torus system $Y'_*$ in $X'$ to the regular neighborhood system 
$N(S(D_*))$ of the 2-sphere system $S(D_*)$ in $\Sigma$, 
there is an orientation-preserving  diffeomorphism 
$\rho''':\Sigma\to \Sigma$   sending 
$N(S(D_*))$ identically  with  
$\rho''' \rho' uS(D'_*)=S(D'_*)$.   
For the orientation-preserving diffeomorphism $h=\rho''' \rho':\Sigma\to \Sigma$, 
the composite embedding $hu: \Sigma^{(0)}\to \Sigma$  preserves the 
O2-sphere basis $(S(D_*),S(D'_*))$ identically. 
This completes the proof of Lemma~3.4. $\square$

\phantom{x}

\noindent{\bf Completion of Proof of  Theorem  1.1.}
Since 
\[(S(D_*),S(D'_*))= (R_*,R'_*), \quad (uS(D_*),uS(D'_*))=(S_*,S'_*), \]
the orientation-preserving diffeomorphism $h$ of $\Sigma$ in Lemma~3.4 sends $(S_*,S'_*)$ to 
 $(R_*,R'_*)$. This completes the proof of Theorem~1.1.
$\square$

\phantom{x}

\noindent{\bf Note 3.6.} The diffeomorphism $h$ in Lemma~3.4 is taken to be isotopic to the identity by using the following lemma ({\it Framed light-bulb isotopy lemma}) based on Gabai's 4D light-bulb theorem, \cite{G} instead of Lemma 3.5.

\phantom{x}

\noindent{\bf Lemma~3.7 (Framed Light-bulb Isotopy Lemma).}
Let $Y_*$ be a system of mutually disjoint 4D solid tori $Y_i\, (i=1,2,\dots,n)$ in
$S^4$. 
Let $D_*\times I$  and $E_*\times I$ be  systems of  mutually disjoint framed disks 
 $D_i\times I, E_i\times I\, (i=1,2,\dots, n)$ in  $\cap Y^c_*$  such that 
 \[(D_*\times I)\cap \partial Y^c_i =(\partial D_i)\times I= (\partial E_i)\times I 
= (E_*\times I)\cap \partial Y^c_i\] 
and $\partial D_i =\partial E_i$  is a  boundary fiber circle  of $Y_i$ for all  $i$. 
If the unions $D_i\cup E_i\, (i=1,2,\dots,n)$ are mutually disjoint, then there is a diffeomorphism 
$h:S^4\to S^4$ which is $Y_*$-relatively  isotopic to the identity such that 
$h(D_*\times I,D_*)=(E_*\times I,E_*)$.

\phantom{x}

\noindent{\bf Proof of Lemma~3.7.}  
As in the proof of Lemma~3.5, assume that a boundary collar system of $D_*$ 
coincides with a boundary collar system of $E_*$. 
First,  show the assertion of the special case $n=1$. 
Let $Y=S^1\times D^3$. Let $D$ and $E$ be proper disks in $Y^c$ 
admitting trivial line bundles 
$D\times I, E\times I$ such that 
\[(\partial D)\times I=(\partial E)\times I \subset \partial Y^c.\]
If the singular  2-sphere $D\cup(-E)$ is not null-homologous in $Y^c$, the disk $D$ is 
$\partial$-relatively homologous to the 2-cycle $E+ m S$ in $(Y^c,\partial Y^c)$ 
for a 2-sphere generator 
$[S]$ of $H_2(Y^c;Z)$ (which is isomorphic to $Z$) 
and some non-zero integer $m$. 
The self-intersection numbers $\mbox{I}([D],[D])$ and $\mbox{I}([E],[E])$  in 
$Y^c$ with  the given framing  of  $\partial D=\partial E$ in $\partial Y^c$ and 
$\mbox{I}([S],[S])$ are all  0.  Thus, 
\[\mbox{I}([D],[D])= \mbox{I}([E],[E])+2m\mbox{I}([S],[E])=0\pm 2m=0\]
and $m=0$. This shows that the singular 2-sphere $D\cup(-E)$ is null-homologous in $Y^c$.  
By Gabai's 4D light-bulb theorem \cite[Theorem~10.4]{G}, there is 
a diffeomorphism $\lambda:Y^c\to Y^c$ such that $\lambda$ is $\partial$-relatively 
isotopic to the identity and $\lambda D=E$, which extends to 
a diffeomorphism $h=\lambda^+:S^4\to S^4$ such that $h$ is $Y$-relatively isotopic to 
the identity and $h D=E$. 
Note that a diffeomorphism of $S^4$ preserves   trivial line bundles on disks 
if the  line bundles on the boundary circles are preserved. This is because  
a sole obstruction that a disk admits a trivial line bundle extending a given line 
bundle on the boundary circle in $S^4$ is that the self-intersection number 
of the disk with the boundary framing given by the line bundle is $0$.  
Thus,  the diffeomorphism $h$ of $S^4$ has $h(D\times I,D)=(E\times I,E)$ and 
the assertion of the special case $n=1$ is shown. 
For the proof in general case, 
let $K=K(k_*)$ be a connected graph in $S^4$ constructed from 
the loop system $k_*$ of the loops $k_i=\partial D_i=\partial E_i,\, (i=1,2,\dots,n)$ by adding mutually disjoint $n-1$ simple arcs $a_j\,(j=1,2,\dots,n-1)$ not meeting 
any interior disk of the disk systems $D_*$ and $E_*$. 
For every $s$ with $1\leq s\leq n$, let $Y_s$ be a regular neighborhood 
of the disk-arc union 
\[k_s\cup_{1\leq i\leq n, i\ne s} E_i\cup_{j=1}^{n-1} a_j\] 
in $S^4$,  which is a 4D solid torus in $S^4$. 
By the proof of the special case $n=1$, 
there is a diffeomorphism $h_1$ of $S^4$ such that $h_1$ is   
$Y_1$-relatively isotopic to the identity  and  $h_1 D_1=E_1$. 
Next, for the disk systems $h_1(D_*)$ and $E_*$, 
there is a diffeomorphism $h_2$ of $S^4$ such that $h_2$ is   
$Y_2$-relatively isotopic to the identity and   
$h_2 h_1 D_1=E_1$ and $h_2 h_1 D_2=E_2$. 
Continuing this process, 
there is a diffeomorphism $h=h_n \dots h_2 h_1$ of $S^4$ such that $h$ is   
$N(K)$-relatively isotopic to the identity for 
a regular neighborhood $N(K)$ of $K$ and 
 $h D_i=E_i\, (i=1,2,\dots,n)$. 
Thus, for a regular neighborhood $N(k_*)$ of the loop system $k_*$ in $N(K)$,
this diffeomorphism $h$ of $S^4$  is  
$N(k_*)$-relatively isotopic to the identity and  $h D_i=E_i\, (i=1,2,\dots,n)$. 
This diffeomorphism  $h$ of $S^4$ is $N(k_*)$-relatively isotopic to the identity for a regular neighborhood $N(k_*)$ of the loop system $k_*$ in $N(K)$ and has $h D_i=E_i\, (i=1,2,\dots,n)$, where $N(k_*)$ can be regarded as the 4D solid torus system $Y_*$ as in the proof of Lemma~3.5.  This completes the proof of Lemma~3.7. $\square$

\phantom{x}

\noindent{\bf 4. Diffeomorphisms of  stable 4-sphere }

Let  $\mbox{Diff}^+(D^4, \mbox{rel}\partial)$ be the orientation-preserving diffeomorphism group of the 4-ball $D^4$ keeping the boundary $\partial D^4$ point-wise fixed. 
An {\it identity-shift} of a  4-manifold  $\Sigma$ is a diffeomorphism 
$\iota:\Sigma\to \Sigma$ obtained from the identity  
$1:\Sigma\to \Sigma$ by replacing the identity on a  4-ball in $\Sigma$ disjoint 
from $F$ with an element of $\mbox{Diff}^+(D^4, \mbox{rel }\partial)$. 
The following result is a main result of this section. 

\phantom{x}

\noindent{\bf Theorem~4.1.} 
Any two homotopic diffeomorphisms of the stable 4-sphere $\Sigma$ are isotopic up to a composition of one diffeomorphism and an identity-shift $\iota$.

\phantom{x}

Because  at  present {\it  it  appears  unknown  whether} 
$\pi_0(\mbox{Diff}^+(D^4, \mbox{rel}\partial))$ {\it is trivial or not}, the identity-shift $\iota$ is needed in Theorem~4.1. However, it is known that any identity-shift $\iota$ is concordant to the identity since $\Gamma_5=0$ by  Kervaire \cite{Ke}). Thus, the following result is a consequence of Theorem~4.1(, whose proof is omitted).

\phantom{x}

\noindent{\bf  Corollary~4.2.} 
Any two homotopic diffeomorphisms of the stable 4-sphere $\Sigma$ are concordant.

\phantom{x}

In Piecewise-Linear Category, every piecewise-linear auto-homeomorphism of the 4-disk keeping the boundary identically is piecewise-linearly $\partial$-relatively isotopic to the  identity, \cite{Hud},  \cite{RS}.  Thus, the following result is a consequence of Theorem 4.1(, whose proof is omitted).

\phantom{x}

\noindent{\bf Corollary~4.3.}  Any two homotopic piecewise-linear auto-homeomorphisms of the stable 4-sphere $\Sigma$ are piecewise-linearly isotopic.

\phantom{x}

The proof  of  Theorem 4.1 is done as follows.

\phantom{x}

\noindent{\bf Proof  of Theorem~4.1.}  Let $f_i\, (i=0,1)$ are homotopic 
diffeomorphisms of $\Sigma=S^4(F)_2$ for a trivial surface-knot $F$ of genus $n$ in $S^4$. Then the composite diffeomorphism $g=f_1^{-1}f_0$ of $\Sigma$ is homotopic to the identity. By Lemmas~3.2, 3.3, 3.4 and Note 3.6, there is a diffeomorphism $h$ of $\Sigma$ isotopic to the identity such that  the composite diffeomorphism $hg$ of $\Sigma$ sends the O2-sphere basis 
$(S(D_*),S(D'_*))$  identically. By the proof of  Lemma~2.1, there is  a proper 1-punctured  trivial surface $P'$ of  $F$ in a 4-ball $A$ such that a region $\Omega(S(D_*),S(D'_*))$ in $\Sigma$ is the double branched covering space $A(P')_2$ of the 4-ball $A$ branched along $P'$ and  the residual region $\Omega^c(S(D_*),S(D'_*))$ is the double branched covering space $A^c(d)_2$ of the 4-ball $A^c=\mbox{cl}(S^4\setminus A)$ branched along the proper trivial disk-knot 
$d=\mbox{cl}(F\setminus P')$ in $A^c$. In this situation, there is a diffeomorphism $h'$ of  $\Sigma$ which is isotopic to the identity such that the composite diffeomorphism $h'hg$ of $\Sigma$ preserve the region $\Omega(S(D_*),S(D'_*))$ identically, which defines a $\partial$-identical diffeomorphism 
$\delta:A^c(d)_2 \to A^c(d)_2$. Since  $A^c(d)_2$ is a 4-ball  and  the lifting $d'$ of the disk $d$  to $A^c(d)_2$  is a trivial disk-knot in $A^c(d)_2$ which is $\partial$-parallel, there is a 
$\partial$-identical diffeomorphism $\delta'$ of $A^c(d)_2$  which is $\partial$-relatively isotopic to the identity such that the composite diffeomorphism $\delta'\delta$ is the identity except for a 
4-ball $U$ in $A^c(d)_2$  disjoint from $d'$. Let $h''$ be the diffeomorphism of  $\Sigma$  defined by 
$\delta'$ and the identity of $\mbox{cl}(\Sigma\setminus A^c(d)_2)$. 
The composite diffeomorphism $h''h'hg$  of  $\Sigma$  preserves $\mbox{cl}(\Sigma\setminus U)$   identically, which is considered as an  identity-shift  $\iota$  of  $\Sigma$. Since  $h''h'h$ is isotopic to  the identity, the diffeomorphism $g=f_1^{-1}f_0$   of  $\Sigma$  is isotopic to $\iota$, and thus the diffeomorphism $f_0$  of $\Sigma$ is isotopic to the composite diffeomorphism 
$ f_1\iota$  of $\Sigma$. This completes  the proof  of Theorem ~4.1. 

\phantom{x}

\noindent{\bf 5. Conclusion}

In an earlier version of this paper (see the research announcement \cite{KA}), the author tried to show Lemma~3.1 with an $\alpha$-equivariant  orientation-preserving diffeomorphism as $f$ (see \cite[Lemma 3.2]{KA}) by a homological argument on an O2-handle basis of a trivial surface-knot $F$ of genus $n$ in $S^4$. However, such a trial is not yet succeeded. The cause of this failure arose from a calculation error on the intersection numbers of O2-handle bases. In fact, the claim 
\cite[Lemma 3.1]{KA} is false, which can be seen by checking the intersection numbers of the 
sphere-bases $(S(D_1),S(D'_1))$ and $(S(E_1),S(E'_1))$ in $\Sigma(1)$, although the $Z_2$-version of \cite[Lemma 3.1]{KA} is true. The claim \cite[Lemma 3.2]{KA} and the related claims in \cite{KA} 
except for \cite[Lemma 3.1]{KA}  are  affirmatively solved  by using Theorem~1.1, Corollary~2.2 and Theorem~4.1 if {\it every sphere-basis} $(S*,S'*)$ {\it of}  $\Sigma$  
{\it is homotopic to the sphere-basis}  $(S(E_*),S(E'_*))$ 
 {\it of an O2-handle basis} 
$(E_*,E'_*)$ 
{\it of} $F$
{\it in} 
$S^4$.  It is hoped that attempts to understand every sphere basis of $\Sigma$ by using  
the O2-handle bases of a trivial surface-knot $F$ in $S^4$  will help simplify the proof of Theorem~1.1.

\phantom{x}

\noindent{\bf Acknowledgments.} 
Since writing the first draft of this paper (in an announcement form) in late November of 2019,
 I stayed at Beijing Jiaotong University, China from December 18, 2019 to January 4, 2020 
while improving the paper. I would like to thank  Liangxia Wan and Rixin Zhang (graduate student) 
for their warm hospitalities. 
This work was partly supported  by JSPS KAKENHI Grant Numbers JP19H01788, JP21H00978  and MEXT Promotion of Distinctive Joint Research Center Program JPMXP0723833165.

\phantom{x}


\begin{thebibliography}{99}

\bibitem{AP} A. Akhmedov and B. D. Park, Geography of simply connected spin 
symplectic 4-manifolds, Math. Res. Lett., 17 (2010), 483-492.

\bibitem{Cerf} J. Cerf, Sur les diff{\'e}omorphismes de la sphere de 
dimension trois ($\Gamma_4=0$), Springer Lecture Notes in Math. 53 (1968). 

\bibitem{Fr} M.  Freedman, The disk theorem for four-dimensional manifolds,  
Proc. Internat. Congr. Math. (Warsaw, Poland) (1983), 647-663. 

\bibitem{FQ} M. H. Freedman and F. Quinn, Topology of 4-manifolds, 
Princeton Univ. Press (1990).

\bibitem{G} D. Gabai, The 4-dimensional light bulb theorem, J. Amer. Math. Soc. 33 (2020), 609-652.


\bibitem{Hat} A. Hatcher, A proof of the Smale conjecture, Diff($S^3$) $\simeq$ O(4),
Ann. of Math.,  117 (1983), 553-607.

\bibitem{HiK} J.  A.  Hillman and  A. Kawauchi, Unknotting orientable surfaces in the 4-sphere, J. Knot Theory Ramifications 4(1995), 213-224.

\bibitem{Hiro}  S. Hirose, On diffeomorphisms over surfaces trivially 
embedded in the 4-sphere, Algebraic and Geometric Topology 2(2002), 791-824. 

\bibitem{HK} F. Hosokawa and A. Kawauchi, Proposals for unknotted 
surfaces in four-space, Osaka J. Math. 16 (1979), 233-248. 

\bibitem{Hud} J. F. P. Hudson, Piecewise-linear topology, Benjamin (1969).

\bibitem{K} A. Kawauchi, Ribbonness of a stable-ribbon surface-link, I. 
A stably trivial surface-link, Topology and its Applications 301(2021), 107522 (16pages). 

\bibitem{KA} A. Kawauchi, Smooth homotopy 4-sphere (research announcement), 2191 Intelligence of Low Dimensional Topology, RIMS Kokyuroku 2191 (July 2021), 1-13. 

\bibitem{KS} A. Kawauchi, Uniqueness of an orthogonal 2-handle pair on a surface-link, 
Contemporary Mathematics (UWP) 4 (2023), 182-188. 

\bibitem{K'} A. Kawauchi,Triviality of a surface-link with meridian-based free fundamental group. Transnational Journal of Mathematical Analysis and Applications 11 (2023), 19-27.


\bibitem{Ke} M. Kervaire, Groups of homotopy spheres: I, Ann. of Math. 77 (1963) 504-537.

\bibitem{Kirby} R. C. Kirby, The topology of 4-manifolds,  Lecture Notes in Mathematics  1374 
(1989),  Berlin: Springer-Verlag. 

\bibitem{RS}  C. Rourke and B. Sanderson, Introduction to piecewise-linear topology, 
Ergebnesse der Math., 69 (1972), Springer-Verlag.

\bibitem{W0} C. T. C. Wall, Diffeomorphisms of 4-manifolds, J. London Math. Soc. 39 (1964),  
131-140.

\bibitem{W} C. T. C. Wall, On simply-connected 4-manifolds, J. London Math. Soc., 39 
(1964), 141-149.



\end{thebibliography}
\end{document}